\documentclass[11pt]{amsart}
\usepackage{amsmath}

\theoremstyle{plain}
\newtheorem{thm}{Theorem}[section]
\newtheorem{lem}[thm]{Lemma}

\theoremstyle{definition}

\def \CPb {\overline{\mathbf{CP}}^{\,2}}
\def \CP {{\mathbf{CP}}^{\,2}}

\def \R {\mathbf{R}}
\def \Z {\mathbf{Z}}
\def \Sig{\Sigma}

\def \la {\langle}
\def \ra {\rangle}

\def \b {\beta}
\def \g {\gamma}

\def \k {\kappa}

\def \o {\omega}
\def \s {\sigma}
\def \t {\tau}

\def \bd {\partial}
\def \x {\times}

\def \sw {\text{SW}}

\def \DD {\Delta}

\begin{document}

\baselineskip.525cm

\title{ Surfaces in 4-Manifolds}
\author[Ronald Fintushel]{Ronald Fintushel}
\address{Department of Mathematics, Michigan State University \newline
\hspace*{.375in}East Lansing, Michigan 48824}
\email{\rm{ronfint@math.msu.edu}}
\thanks{The first author was partially supported NSF Grant DMS9401032 and
the second
author by NSF Grant DMS9626330}
\author[Ronald J. Stern]{Ronald J. Stern}
\address{Department of Mathematics, University of California \newline
\hspace*{.375in}Irvine,  California 92697}
\email{\rm{rstern@math.uci.edu}}
\date\today

\begin{abstract} In this paper we introduce  a  technique, called {\em rim surgery}, which can
change a smooth embedding of an
orientable surface $\Sig$ of positive genus and nonnegative
self-intersection in a smooth
$4$-manifold $X$ while leaving the topological embedding unchanged. This is accomplished by replacing
the tubular neighborhood of a
particular nullhomologous torus in $X$ with $S^1\x E(K)$, where
$E(K)$ is the exterior of a knot $K\subset S^3$. The smooth change can be detected
easily for certain
pairs $(X,\Sig)$ called {\em SW-pairs}. For example,  $(X,\Sig)$  is an
SW-pair if $\Sig$ is a
symplectically and primitively embedded surface with positive genus and
nonnegative
self-intersection in a simply connected symplectic 4-manifold $X$. We prove
the following
theorem:
\smallskip

\noindent{\bf Theorem.} {\em Consider any SW-pair $(X,\Sig)$. For each knot
$K\subset S^3$
there is a surface $\Sig_K\subset X$ such that the pairs $(X,\Sig_K)$  and
$(X,\Sig)$ are homeomorphic. However, if $K_1$ and $K_2$  are two knots for
which there
is a diffeomorphism of pairs $(X, \Sig_{K_1})\to (X,\Sig_{K_2})$, then
their Alexander
polynomials are equal:
$\DD_{K_1}(t)=\DD_{K_2}(t)$.}
\end{abstract}

\maketitle
\section{Introduction\label{Intro}} We say that a surface $\Sig$ is  {\it
primitively
embedded}  in a simply connected smooth 4-manifold
$X$ if $\Sig$ is smoothly embedded with $\pi_1(X\setminus \Sig)=0$. In
particular, by
Alexander  duality, $\Sig$ must represent a primitive homology class $[\Sig] \in
H_2(X;\Z)$. In general, any smoothly embedded (connected) surface $S$ in a
simply
connected smooth 4-manifold $X$  with $[S]\ne 0$ has the property that the
surface
 $\Sig$ which represents the homology class $[S]-[E]$ in $X\#\CPb$ and which is
obtained by tubing
together the surface
$S$ with the exceptional sphere $E$ of $\CPb$ is primitively embedded (since
the surface $\Sig$ transversally intersects the sphere $E$ in one point).

Given  a primitively embedded positive genus surface $\Sig$ in $X$, in the
first part of
this paper we shall construct for each knot
$K$ in $S^3$ a smoothly embedded surface
$\Sig_K$ in $X$ which is {\it $\Sig$-compatible}; i.e. $[\Sig]=[\Sig_K]$
and there is a
homeomorphism $(X, \Sig)\to (X,\Sig_K)$. This construction will have two
properties.
The first is  that
$(X,\Sig_{\text{unknot}})=(X,\Sig)$. The main result of this paper is the second
property: under suitable hypotheses on the pair $(X,\Sigma)$,  if
$K_1$ and $K_2$ are two knots in $S^3$ and if there is a diffeomorphism
$(X, \Sig_{K_1})\to (X,\Sig_{K_2})$, then
$K_1$ and $K_2$ have the same symmetric Alexander polynomial, i.e.
$\DD_{K_1}(t)=\DD_{K_2}(t)$. As a special case we  show:

\begin{thm}\label{sympthm} Let $X$ be a simply connected symplectic
4-manifold and
$\Sig$ a symplectically and primitively embedded surface with positive genus and
nonnegative self-intersection. If
$K_1$ and $K_2$ are knots in $S^3$
 and if there is a diffeomorphism of pairs
$(X, \Sig_{K_1})\to (X,\Sig_{K_2})$, then  $\DD_{K_1}(t)=\DD_{K_2}(t)$.
Furthermore,
if $\DD_K(t)\ne 1$, then
$\Sig_K$ is not smoothly ambient isotopic to a symplectic submanifold of $X$.
\end{thm}

For example, Theorem~\ref{sympthm} applies to the $K3$ surface where $\Sig$ is a
generic elliptic fiber. It also applies to surfaces of the form $S - E$ in
$\CP \# \CPb$,
where $S$ is any positive genus symplectically
embedded surface in $\CP$.

The outline of this paper is as follows. In
\S 2 we shall construct the surfaces
$\Sig_K$ with $[\Sig_K]=[\Sig]$ and  show that if
$\pi_1(X)=\pi_1(X\setminus \Sig)=0$,
there is a homeomorphism of $(X, \Sig)$ with $(X,\Sig_K)$, i.e. $\Sig_K$ is
$\Sig$-compatible.  We give two descriptions of $\Sig_K$.  One is explicit,
while the other describes how to obtain $\Sig_K$ by removing a tubular
neighborhood
$T^2\x D^2$ of a homologically  trivial torus in a tubular neighborhood  of
$\Sig$ and
replacing it with $S^1 \x E(K)$, where $E(K)$ is the exterior  of the  knot
$K$ in $S^3$.
This is reminiscent of our  construction in \cite{FS} where we performed
the same
operation on homologically essential tori.
There, the Alexander polynomial $\DD_K(t)$ of $K$ detected  a change in the
diffeomorphism type of the ambient manifold $X$. Here, we shall show that
$\DD_K(t)$
detects a change in the diffeomorphism type of the  embedding of $\Sig$ in
$X$.

If the self-intersection of $\Sig$ is $n\ge 0$, then in $X_n=X\#n\CPb$
consider the
surface
$\Sig_n=\Sig-\sum_{j=1}^n E_j$ (resp.
$\Sig_{n,K}=\Sig_K-\sum_{j=1}^n E_j$) obtained from
$\Sig$  (resp. $\Sig_K$) by tubing together with the exceptional spheres $E_j$,
$j=1,\dots,n$, of the copies of $\CPb$ in
$X_n$. If there is a diffeomorphism
$H:(X, \Sig_{K_1})\to (X,\Sig_{K_2})$, then there is a diffeomorphism
$H_n:(X_n, \Sig_{n,K_1})\to (X_n,\Sig_{n,K_2})$.  For each genus $g\ge1$
we construct in
\S3 a standard pair $(Y_g, S_g)$, with the properties that $Y_g$ is a
Kahler surface,
$S_g$ is a primitively embedded genus $g$ Riemann surface in $Y_g$, and
the torus
used to construct $S_{g,K}=(S_g)_K$ is contained in a cusp neighborhood.
Then in \S4 we will study {\it SW-pairs}, i.e.  pairs $(X,\Sig)$ where
$X$ is a smooth simply connected  4-manifold, $\Sig$ is a
primitively embedded genus $g$ surface with self-intersection $n\ge 0$, and
the fiber sum of
$X_n$ and $Y_g$ along the surfaces $\Sig_n$ and $S_g$ has a nontrivial
Seiberg-Witten
invariant $\sw_{X_n\#_{\Sig_n=S_g} Y_g}\ne 0$.  The point here is that the
nullhomologous
torus used to construct the surface $\Sig_K$ in $X$ still resides in
$X_n\#_{\Sig_n=S_g} Y_g$ and is now homologically essential and is
contained in a
cusp neighborhood. It will also follow that if
$X$ is a symplectic 4-manifold and
$\Sig$ is a symplectically and primitively embedded surface with nonnegative
self-intersection, then $(X,\Sig)$ is a SW-pair.

In \S 5 we use in a straightforward fashion the results of \cite{FS}
to show that the Alexander polynomial of
$K$ distinguishes the
$\Sig_K$ for SW-pairs, and we complete the proof our main theorem:

\begin{thm}\label{mainthm} Consider any $\sw$-pair $(X,\Sig)$. If
$K_1$ and $K_2$ are two knots in $S^3$ and if there is a diffeomorphism of pairs
$(X, \Sig_{K_1})\to (X,\Sig_{K_2})$, then  $\DD_{K_1}(t)=\DD_{K_2}(t)$.
\end{thm}
\noindent Finally, in \S6 we complete the proof of Theorem~\ref{sympthm}
by showing that
in the case  that $\Sig$ is sympletically embedded in $X$ and  $\DD_K(t)\ne
1$, then
$\Sig_K$ is not smoothly ambient isotopic to a symplectic submanifold of
$X$.

We  conclude this introduction with  two conjectures. The first conjecture
is that, under
the hypothesis of Theorem~\ref{mainthm},  there is a diffeomorphism
$(X, \Sig_{K_1})\to (X,\Sig_{K_2})$ if and only if the knots $K_1$ and $K_2$ are
isotopic.  In particular, this conjecture would
imply that the study of the equivalence clases of $\Sig$-compatible
surfaces under
diffeomorphism is at least as complicated as classical knot theory. The
second conjecture
is a finiteness conjecture: given a symplectic 4-manifold
$X$ and a symplectic submanifold $\Sig$, we conjecture that there are only
finitely many
distinct smooth isotopy classes of symplectic submanifolds $\Sig^\prime$
which are
topologically isotopic to $\Sig$.

\section{The construction of $\Sig_K$}

Let $X$ be a smooth
$4$-manifold which contains a smoothly embedded surface $\Sig$ with genus $g>0$.
Then there is a diffeomorphism
\[ h:\Sig \to T^2\#\cdots\# T^2 =
(T^2\setminus D^2)\cup (T^2\setminus(D^2 \amalg D^2))\cup \cdots\cup
(T^2\setminus
D^2).\] Let $C\subset \Sig$ be a curve whose image under $h$ is the curve
$S^1\x\{\text{pt} \} \subset T^2\setminus D^2 = (S^1\x S^1)\setminus D^2$
in the first
$T^2\setminus D^2$ summand of
$h(\Sig)$.  Keep in mind that, since there are many such diffeomorphisms
$h$, there are many such curves $C$. Given a knot $K$ in $S^3$
we shall give two different constructions of a surface
$\Sig_{K,C}$. The first is an explicit construction, while the  second
shows how to
obtain $\Sig_{K,C}$ by what we call a  {\it{rim surgery}},  a surgical
operation on a particular homologically trivial torus in a
neighborhood of $\Sig$. It is this second  construction that will allow us
to compute
appropriate invariants to distinguish the surfaces $\Sig_{K,C}$.
\subsection{An explicit description of $\Sig_{K,C}$}
 Viewing $S^1$ as the union of two arcs $A_1$ and $A_2$, we have
\begin{eqnarray*}T^2\setminus D^2 &=& (S^1\x S^1)\setminus D^2\\&=&((A_1\cup
A_2)\x(A_1\cup A_2))
\setminus (A_1 \times A_1)
\\&=&(A_2\x S^1)\cup (A_1 \x A_2)
\end{eqnarray*} with $h(C)=A_2\x\{{\text{pt}}\}\cup A_1\x\{{\text{pt}}\}$.
Now the normal
bundle of $\Sig$ in $X$ when restricted to
$T^2\setminus D^2\subset \Sigma$ is trivial, hence it is diffeomorphic to
\[ ((A_2\x S^1)\cup (A_1 \x A_2))\x D^2=((A_2\x D^2)\x S^1 )\cup ((A_1\x
D^2) \x A_2)).\]
Furthermore, under this diffeomorphism, the inclusion
\[ (T^2\setminus D^2)\x \{0\} \subset  (T^2\setminus D^2)\x D^2 \]
becomes
\[ (A_2\x \{0\})\x S^1)\cup ((A_1 \x \{0\})\x A_2) \subset ((A_2\x D^2)\x
S^1)\cup
((A_1\x D^2)\x A_2).\]
Now tie a knot $K$ in  the arc  $(A_2\x \{0\}) \subset (A_2\x D^2)$ to
 obtain a knotted  arc $A_K$ and  to obtain a new punctured torus
\begin{eqnarray*}   T_K\setminus D^2  &= &(A_K\x S^1)\cup ((A_1 \x \{0\})\x A_2)
\\&\subset&((A_2\x D^2)\x S^1)\cup ((A_1 \x D^2)\x A_2)
\end{eqnarray*} with \[\partial (T_K\setminus D^2)= \partial (T\setminus
D^2).\] Then let
\[\Sig_{K,C}
 = (T_K\setminus D^2)\cup(T^2\setminus(D^2\amalg D^2))\cup\cdots\cup
(T^2\setminus D^2)
\subset N(\Sigma)\subset X.\]

\subsection{A description of $\Sig_{K,C}$ via rim surgery} Keeping  the
notation above,  we first recall how, via a 3-manifold surgery, we can tie
a knot $K$ in
the arc  $(A_2\x \{0\}) \subset (A_2\x D^2)$.
In short, we just remove a small tubular neighborhood in $A_2\x D^2$ of a
pushed-in
copy $\g$ of the meridional circle
$\{0\}\x S^1 \subset A_2\x D^2$ and sew in the exterior of the knot $K$ in $S^3$
so that the meridian of $K$ is identified with $\g$. This has the effect of
tying
a knot in the arc $A_2\x\{0\}\subset A_2\x D^2$. More specifically,
consider the standard
embedding of the solid torus
$A=(A_1\cup A_2)\x D^2=S^1\x D^2$ in $S^3$ with complementary solid torus
$B=D^2\x S^1$
with core
$C^\prime = \{0\} \x S^1\subset D^2 \x S^1$. In   $A\setminus C=(S^1 \x
D^2)\setminus
C=S^1\x S^1\x (0,1]=(A_1 \cup A_2)\x S^1\x (0,1]$, consider the circle
$\g= \{t\} \x  S^1 \x \{\frac{1}{2}\}$, with $t \in A_2$, and with tubular
neighborhood
$N(\g)\subset A\setminus C$. The curve $\g$ is isotopic in $S^3\setminus C$
to the core
$C^\prime$ of $B$.  We denote by $\g^\prime $ the curve
$\g$ pushed off into
$\partial N(\g)$ so that the linking number in $S^3$ of $\g$ and
$\g^\prime$ is zero.  For later reference, note that
$D= (A\setminus N(\g))\cup B$ is again diffeomorphic to a solid torus.
(It is the exterior of the unknot
$\g \subset A \subset S^3$.) The core of $D$ is isotopic (in $D$) to $C$.

Let $M_K$ be the 3-manifold obtained by performing $0$-framed surgery on
$K$. Then the
meridian $m$ of $K$ is a circle in $M_K$ and has a  canonical framing in
$M_K$; we denote a tubular neighborhood of $m$ in $M_K$ by $m\x D^2$. Let
$S_K$ denote the 3-manifold
\[ S_K=(A\setminus N(\g)) \cup (M_K\setminus (m\x D^2)).\] The two pieces
are glued
together so as to identify $\g^\prime$ with $m$. In other words, we remove
$N(\g)$ and sew
in the exterior
$E(K)$ of the knot $K$ in $S^3$.  Note that the core $C$ of the  solid
torus $A$ is
untouched by this  operation, so
$C \subset S_K$. Also, the boundary $\partial A$ of $A$ and the set
$G=A_1\x D^2 \subset (A_1\cup A_2)\x D^2 \subset A$ remain untouched and
thus can be
viewed as subsets of $S_K$.

\begin{lem}\label{knot} There is a diffeomorphism $h: S_K  \to A$  which is
the identity
on $G$ and on the boundary. Furthermore,
$h(C)$ is the knotted core $K \subset A$.
\end{lem}
\begin{proof} In $S^3=A \cup B$, the above operation replaces a tubular
neighborhood of
the unknot $\g \subset A \subset S^3$ with the  exterior $E(K)$ of the knot
$K$ in $S^3$.
Thus there is a diffeomorphism
$h: E(K)\cup D \to A\cup B=S^3$ sending the core circle of $D$ to the knot
$K$.
Now $C'\subset B\subset E(K)\cup D$ is unknotted, since in $D$, the curve
$C'$ is isotopic to $\g'$, which bounds a disk.
Thus  $S_K$, which is the complement of a tubular neighborhood of $C'$,
is an unknotted solid torus in $S^3=E(K)\cup D$. Furthermore, as we have
noted above,
$C$ is isotopic to the core of $D$; so
$C\subset S_K$ is the knot $K$. Thus there is a diffeomorphism $h: S_K  \to
S^1\x D^2$ which is the identity on the boundary. After an isotopy rel boundary
we can arrange that $h(G)=G$.
\end{proof}

To obtain $\Sig_{K,C}$ we cross everything with $S^1$; i.e.  remove the
neighborhood
 $N(\g)\x S^1 \subset  (A_2\x D^2) \x S^1 \subset N(\Sig)$ of the
(nullhomologous) torus
$\g\x S^1 \subset (A_2\x D^2) \x S^1 \subset N(\Sig)$ and sew in $E(K)\x
S^1$ as above on
the $E(K)$ factor and the identity on the $S^1$ factor.
 We refer to this as  a {\it rim surgery} on $\Sigma$. Notice that this
construction does
not change the ambient manifold $X$. Except where it is absolutely
necessary to keep
track of the curve $C$, we shall suppress it from our notation and
abbreviate $\Sig_{K,C}$
as $\Sig_K$.

\subsection{The complement of $\Sig_K$} From the construction,
 it is clear that if the complement of $\Sig$ in $X$ is simply connected,
then so is the
complement of $\Sig_K$ in $X$, since  the meridian of the knot (which is
identified with
the boundary of the normal fiber to $\Sig$) normally generates the
fundamental group of
the exterior of $K$. Now there is a map $f: E(K)\to B\cong D^2\x S^1$ which
induces isomorphisms on homology and restricts to a homeomorphism $\bd
E(K)\to\bd B$ taking the class of a meridian to $[\{ {\text{pt}} \}\x S^1]$
and the class
of a longitude to $[\bd D^2\x\{ {\text{pt}}\} ]$. The map $f\x {\text
{id}}_{S^1}$ on
$E(K)\x S^1$ extends via the identity to a homotopy
equivalence $X\setminus N(\Sig_K)\to X\setminus N(\Sig)$ which restricts to a
homeomorphism $\bd N(\Sig_K)\to \bd N(\Sig)$. Then topological surgery
\cite{Fr, B} guarantees the existence of a homeomorphism $h:(X,\Sig)\to
(X,\Sig_K)$.

If $\pi_1(X\setminus\Sig)\ne 0$, it is not clear when $X\setminus \Sig_K$
is homeomorphic
(or even homotopy equivalent) to $X\setminus \Sig_K$. We avoid such issues
in this paper
and only deal with the case where
$\pi_1(X\setminus\Sig)=0$.  However, as already noted; the surface $\Sig -
E$ in $X\#
\CPb$  obtained by tubing together  the surface $\Sig$ with the exceptional
sphere $E$ of
$\CPb$ is primitively embedded; so there is a homeomorphism
 $h:(X\# \CPb, \Sig-E)\to (X\# \CPb,\Sig_K-E)$. In summary:
\begin{thm} Let $X$ be a simply connected smooth 4-manifold with a
 primitively embedded surface $\Sig$. Then for each knot $K$ in
$S^3$, the above construction produces a $\Sig$-compatible surface $\Sig_K$.
\end{thm}

\section{The standard pair $(Y_g,S_g)$ }

Let $g>0$. In this section we shall construct a simply connected smooth
4-manifold $Y_g$
and a primitive embedding of $S_g$, the surface of genus $g$, in $Y_g$ such
that the
torus used in the previous section to construct the $S_g$-compatible embedding
$(S_g)_K=S_{g,K}$
is contained in a cusp neighborhood.

To this end, consider the $(2,2g+1)$-torus knot $T(2,2g+1)$. It is a
fibered knot
whose fiber is a punctured genus $g$ surface and whose monodromy $t'$ is
periodic of
order $4g+2$. If we attach a 2-handle to $\bd B^4$ along $T(2,2g+1)$ with
framing
$0$, we obtain a manifold $C(g)$ which fibers over the 2-disk with generic
fiber a
Riemann surface $S_g$ of genus $g$ and whose monodromy map $t$, induced
from $t'$,
is a periodic holomorphic map $t:S_g\to S_g$ of order $4g+2$. The singular fiber
is the topologically (non-locally flatly) embedded sphere obtained from the
cone in $B^4$ on
the torus knot $T(2,2g+1)$ union the core of the 2-handle. Now consider the
fibration
over the punctured 2-sphere obtained from gluing together
 $4g+2$ such neighborhoods $C(g)$ along a neighborhood of a fiber in the
boundary of
$C(g)$. This is a complex
surface, and the
monodromy is trivial around a loop which contains in its interior the
images of all
the singular fibers. Thus we may compactify this manifold to obtain a complex
surface $Y_g$ which is holomorphically fibered over $S^2$. For example,
$Y_1$ is the
rational elliptic surface $\CP\#9\CPb$. In fact, $Y_g$ is just the Milnor
fiber of
the Brieskorn singularity $\Sigma(2,2g+1,4g+1)$ union a generalized nucleus
consisting of the
4-manifold obtained as the trace of the 0-framed surgery on $T(2,2g+1)$ and
a $-1$ surgery on a
meridian \cite{Fuller}. The fibration
$\pi: Y_g\to S^2$ has a holomorphic section which is a sphere $\Lambda$ of
self-intersection
$-1$  (the sphere obtained by the $-1$-surgery above (cf. \cite{Fuller}).
This proves
that
$\pi_1(Y_g\setminus S_g)=0$; so $S_g$ is a primitively embedded surface with
self-intersection $0$.

Let $T$ denote the torus in $S_g\x D^2$ on which we perform a rim surgery
in order to obtain the surface $S_{g,K}$. We wish to see that $T$ lies in a cusp
neighborhood.
A cusp neighborhood is nothing more than the regular neighborhood of a
torus together with two
vanishing cyles,
one for each generating circle in the torus. The torus $T$ has the form
$T=\g\x\t$ where
$\t$ is a closed curve on $S_g$ and $\g = \{{\text {pt}}\}\x
(\{\frac12\}\x\bd D^2)$. The
curve $\t$ is one of the generating circles for $H_1(S_g;\Z)$ with a
dual circle $\s$. The curve $\g$ spans a $-1$-disk contained in $\Lambda$.
The curve
$\t$ degenerates to a point on the singular fiber in $C(g)$. Thus we
see both required vanishing cycles.

\section{SW-pairs}

Recall that the Seiberg-Witten invariant  $\sw_X$ of a smooth closed
oriented $4$-manifold
$X$ with $b ^+>1$ is an integer valued function which is defined on the set of
$spin^{\, c}$ structures over $X$, (cf. \cite{W}). In case $H_1(X;\Z)$ has no
2-torsion, there is a natural identification of the
$spin^{\, c}$ structures of
$X$ with the characteristic elements of $H^2(X;\Z)$. In this case we view the
Seiberg-Witten invariant as
\[ \sw_X: \lbrace k\in H^2(X,\Z)|k\equiv w_2(TX)\pmod2)\rbrace
\rightarrow \Z. \] The Seiberg-Witten invariant $\sw_X$ is a smooth
invariant whose sign
depends on an orientation of
$H^0(X;\R)\otimes\det H_+^2(X;\R)\otimes \det H^1(X;\R)$. If $\sw_X(\b)\neq
0$, then we
call $\b$ a {\it{basic class}} of $X$. It is a fundamental fact that the
set of basic
classes is finite.  If $\b$ is a basic class, then so is $-\b$ with
\[\sw_X(-\b)=(-1)^{(\text{e}+\text{sign})(X)/4}\,\sw_X(\b)\] where
$\text{e}(X)$ is the Euler number and $\text{sign}(X)$ is the signature of $X$.

As in \cite{FS} we need to view the Seiberg-Witten invariant
 as a Laurent polynomial. To do this,  let
$\{\pm \b_1,\dots,\pm \b_n\}$ be the set of nonzero basic classes for $X$.
We my then
view the Seiberg-Witten invariant of $X$ as the `symmetric' Laurent polynomial
\[\sw_X = b_0+\sum_{j=1}^n
b_j(t_j+(-1)^{(\text{e}+\text{sign})(X)/4}\,t_j^{-1})\]
where  $b_0=\sw_X(0)$, $b_j=\sw_X(\b_j)$ and $t_j=\exp(\b_j)$.

Now let $\Sig$ be genus $g>0$ primitively embedded surface in the simply
connected
4-manifold $X$. If the self-intersection of $\Sig$ is $n\ge
0$, then in
$X_n=X\#n\CPb$, consider the surface
$\Sig_n=\Sig-\sum_{j=1}^n E_j$ (resp.
$\Sig_{n,K}=\Sig_K-\sum_{j=1}^n E_j$) obtained from
$\Sig$  (resp. $\Sig_K$) by tubing together with the  exceptional spheres
$E_j$, $j=1,\dots,n$,  of the $\CPb$ in $X_n$. Note that the fiber sum
$X_n\#_{\Sig_n=S_g} Y_g$ of $X_n$ and $Y_g$ along the surfaces $\Sig_n$ and
$S_g$ has $b^+>1$.
An {\it SW-pair}\/ is such a pair $(X,\Sig)$ which satisfies the property
that the
Seiberg-Witten invariant $\sw_{X_n\#_{\Sig_n=S_g} Y_g}\ne 0$.

As we have pointed out earlier, there are several curves $C$ that can be
used to construct
the  surfaces $\Sig_{K,C}$, and there are potentially several different
fiber sums that
can be performed in the construction of
$X_n\#_{\Sig_n=S_g} Y_g$. We pin down our choice of $C$ by declaring it to
be the
image of the curve $\s$ from \S 3 under the diffeomorphism used in the
construction of the
fiber sum. A simple Mayer-Vietoris argument shows that in
$X_n\#_{\Sig_n=S_g} Y_g$
the rim torus (equivalently
$\g\x\t$) becomes homologically essential and is contained in a cusp
neighborhood. Thus our
results from \cite{FS} apply.

\section{SW-pairs and the Alexander polynomial \label{proof}}

We are now in a position to prove our main theorem:
\smallskip

\noindent{\bf{Theorem 1.2.}} {\em Consider any $\sw$-pair $(X,\Sig)$. If
$K_1$ and $K_2$ are two knots in $S^3$ and if there is a diffeomorphism of pairs
$(X,\Sig_{K_1})\to (X,\Sig_{K_2})$, then  $\DD_{K_1}(t)=\DD_{K_2}(t)$.}

\begin{proof} With notation as above, we have a diffeomorphism
$(X_n,\Sig_{n,K_1})\to (X_n,\Sig_{n,K_2})$. Then there is a diffeomorphism
\[ Z_1=X_n\#_{\Sig_{n,K_1}=S_g}Y_g\to Z_2=X_n\#_{\Sig_{n,K_2}=S_g}Y_g. \]
It follows from \cite{FS} that
$\sw_{Z_i}=\sw_{X_n\#_{\Sig_n=S_g} Y_g}\cdot\DD_{K_i}(t)$
for $t=\exp(2[T])$, where $T$ denotes the rim torus. Since $(X,\Sig)$ is a
$\sw$-pair, and
since $[T]\ne 0$ in $H_2(Z_i;\Z)$ we must have $\DD_{K_1}(t)=\DD_{K_2}(t)$.
\end{proof}

\section{Rim surgery on symplectically embedded surfaces}

We conclude with a proof of our claim of the introduction.
\smallskip

\noindent{\bf{Theorem 1.1.}} {\em Let $X$ be a simply connected symplectic
4-manifold and
$\Sig$ a symplectically and primitively embedded surface with positive genus and
nonnegative self-intersection. If $K_1$ and $K_2$ are knots in $S^3$
and if there is a diffeomorphism of pairs
$(X,\Sig_{K_1})\to (X,\Sig_{K_2})$, then  $\DD_{K_1}(t)=\DD_{K_2}(t)$.
Furthermore,
if $\DD_K(t)\ne 1$, then
$\Sig_K$ is not smoothly ambient isotopic to a symplectic submanifold of $X$.}

\begin{proof} Since $\Sig$ and $S_g$ are symplectic submanifolds of $X$ and
$Y_g$, the
fiber sum  $X_n\#_{\Sig_n=S_g}Y_g$ is also a symplectic manifold
\cite{Gompf}. Thus
$\sw_{X_n\#_{\Sig_n=S_g}Y_g}\ne 0$ \cite{T1}; so $(X,\Sig)$ forms an
SW-pair. This proves
the first statement of the theorem.

Next, suppose that $\Sig_K$ is smoothly ambient
isotopic to a symplectic submanifold $\Sig'$ of $X$. This isotopy carries
the rim torus
$T$ to a rim torus $T'$ of $\Sig'$. We have
\begin{equation}
\sw_{X_n\#_{\Sig'_n=S_g}Y_g}=\sw_{X_n\#_{\Sig_{n,K}=S_g}Y_g}=
\sw_{X_n\#_{\Sig_n=S_g}Y_g}\cdot\DD_K(t)
\label{invt}\end{equation}
with $t=\exp(2[T'])$ when this expression is viewed as
$\sw_{X_n\#_{\Sig'_n=S_g}Y_g}$.
As above, $[T']\ne 0$ in $H_2(X_n\#_{\Sig'_n=S_g}Y_g;\Z)$.

Symplectic forms $\o_X$ on $X_n$ (with respect to which $\Sig'_n$ is
symplectic) and
$\o_Y$ on $Y_g$ induce a symplectic form $\o$ on the symplectic fiber sum
$X_n\#_{\Sig'_n=S_g}Y_g$ which agrees with $\o_X$ and
$\o_Y$ away from the region where the manifolds are glued together. In
particular, since
$T'$ is nullhomologous in $X_n$, we have $\la\o,T'\ra = \la\o_X,T'\ra=0$.
Now \eqref{invt}
implies that the basic classes of $X_n\#_{\Sig'_n=S_g}Y_g$  are exactly the
classes
$b+2mT'$ where $b$ is a basic class of $X_n\#_{\Sig_n=S_g}Y_g$ and $t^m$ has a
nonzero coefficient in $\DD_K(t)$.  Thus the basic classes of
$X_n\#_{\Sig'_n=S_g}Y_g$ can be grouped into collections
${\mathcal{C}}_{b}=\{ b+2mT'\}$, and if
$\DD_K(t)\ne 1$ then each ${\mathcal{C}}_{b}$ contains
more than one basic class. Note, however, that
$\la\o,b+2mT'\ra=\la\o,b\ra$.
Now Taubes has shown
\cite{T2} that the canonical class $\k$ of a symplectic manifold with
$b^+>1$ is the
basic class which
is characterized by the condition
$\la\o,\k\ra >\la\o,b'\ra$ for any other basic class $b'$. But this is
impossible for
$X_n\#_{\Sig_n=S_g}Y_g$ since each ${\mathcal{C}}_{b}$ contains more than
one class.
\end{proof}

\end{document}